\documentstyle[11pt,twoside]{article}
\include{graphicx}
\title{A note on an integral of Dixit, Roy and Zaharescu}
\author{\sc R. B.\ Paris \\
{\em Division of Computing and Mathematics}, \\
{\em Abertay University, Dundee DD1 1HG, UK}
}
\begin{document}
\def\f#1#2{\mbox{${\textstyle \frac{#1}{#2}}$}}
\def\dfrac#1#2{\displaystyle{\frac{#1}{#2}}}
\def\boldal{\mbox{\boldmath $\alpha$}}
{\newcommand{\Sgoth}{S\;\!\!\!\!\!/}
\newcommand{\bee}{\begin{equation}}
\newcommand{\ee}{\end{equation}}
\newcommand{\lam}{\lambda}
\newcommand{\ka}{\kappa}
\newcommand{\al}{\alpha}
\newcommand{\fr}{\frac{1}{2}}
\newcommand{\fs}{\f{1}{2}}
\newcommand{\g}{\Gamma}
\newcommand{\br}{\biggr}
\newcommand{\bl}{\biggl}
\newcommand{\ra}{\rightarrow}
\newcommand{\mbint}{\frac{1}{2\pi i}\int_{c-\infty i}^{c+\infty i}}
\newcommand{\mbcint}{\frac{1}{2\pi i}\int_C}
\newcommand{\mboint}{\frac{1}{2\pi i}\int_{-\infty i}^{\infty i}}
\newcommand{\gtwid}{\raisebox{-.8ex}{\mbox{$\stackrel{\textstyle >}{\sim}$}}}
\newcommand{\ltwid}{\raisebox{-.8ex}{\mbox{$\stackrel{\textstyle <}{\sim}$}}}
\renewcommand{\topfraction}{0.9}
\renewcommand{\bottomfraction}{0.9}
\renewcommand{\textfraction}{0.05}
\newcommand{\mcol}{\multicolumn}
\date{}
\maketitle
\pagestyle{myheadings}
\markboth{\hfill \sc R. B.\ Paris  \hfill}
{\hfill \sc  \hfill}
\begin{abstract}
In a recent paper, Dixit {\it et al.\/} [Acta Arith. {\bf 177} (2017) 1--37] posed two open questions whether the integral
\[{\hat J}_{k}(\alpha)=\int_0^\infty\frac{xe^{-\alpha x^2}}{e^{2\pi x}-1}\,{}_1F_1(-k,\f{3}{2};2\alpha x^2)\,dx\]
for $\alpha>0$ could be evaluated in closed form when $k$ is a positive even and odd integer. We establish that ${\hat J}_{k}(\alpha)$ can be expressed in terms of a Gauss hypergeometric function and a ratio of two gamma functions, together with a remainder expressed as an integral. An upper bound on the remainder term is obtained, which is shown to be exponentially small as $k$ becomes large when $a=O(1)$.
\vspace{0.4cm}

\vspace{0.3cm}

\end{abstract}

\vspace{0.3cm}

\noindent $\,$\hrulefill $\,$

\vspace{0.2cm}

\begin{center}
{\bf 1. \  Introduction}
\end{center}
\setcounter{section}{1}
\setcounter{equation}{0}
\renewcommand{\theequation}{\arabic{section}.\arabic{equation}}
In the first of his letters to Hardy \cite{R1}, Ramanujan gave the formula
\[I(\alpha):=\alpha^{-1/4}\bl(1+4\alpha\int_0^\infty\frac{x e^{-\alpha x^2}}{e^{2\pi x}-1}\,dx\br)=
\beta^{-1/4}\bl(1+4\beta\int_0^\infty\frac{x e^{-\beta x^2}}{e^{2\pi x}-1}\,dx\br),\]
where $\alpha\beta=\pi^2$, and in \cite{R2} obtained the approximate evaluation
\bee\label{e11}
I(\alpha)\simeq\bl(\frac{1}{\alpha}+\frac{1}{\beta}+\frac{2}{3}\br)^{1/4}.
\ee
This approximation is found to be good for small and large values of $\alpha$. A proof of this result was given in \cite{BE}, where the asymptotic expansion
\[I(\alpha)\sim \frac{1}{\alpha^{1/4}}+\frac{\alpha^{3/4}}{6}-\frac{\alpha^{7/4}}{60}+\cdots\qquad (\alpha\to 0)\]
was obtained.

In a recent paper, Dixit, Roy and Zaharescu \cite{DRZ} established an analogous formula for the integral
\[{\hat J}_k(\alpha):=\int_0^\infty\frac{xe^{-\alpha x^2}}{e^{2\pi x}-1}\,{}_1F_1(-k;\f{3}{2};2\alpha x^2)\,dx,\]
where ${}_1F_1$ denotes the confluent hypergeometric function and $k$ is a positive integer. They showed that  \cite[(1.25), (1.27)]{DRZ}
\[\alpha^{-1/4}{}_2F_1(-2k,1;\f{3}{2};2)+4\alpha^{3/4} {\hat J}_{2k}(\alpha)
=\beta^{-1/4}{}_2F_1(-2k,1;\f{3}{2};2)+4\beta^{3/4} {\hat J}_{2k}(\beta)\]
and
\bee\label{e12}
\alpha^{-1/4}{}_2F_1(-2k\!-\!1,1;\f{3}{2};2)+4\alpha^{3/4} {\hat J}_{2k+1}(\alpha)
=-\beta^{-1/4}{}_2F_1(-2k\!-\!1,1;\f{3}{2};2)-4\beta^{3/4} {\hat J}_{2k+1}(\beta)
\ee
when $\alpha\beta=\pi^2$,
where ${}_2F_1$ denotes the Gauss hypergeometric function. 
In the particular case $\alpha=\beta=\pi$, (\ref{e12}) yields the beautiful exact evaluation \cite[Cor. 1.8]{DRZ}
\bee\label{e14}
{\hat J}_{2k+1}(\pi):=\int_0^\infty\frac{xe^{-\pi x^2}}{e^{2\pi x}-1}\,{}_1F_1(-2k\!-\!1;\f{3}{2};2\pi x^2)\,dx=
-\frac{1}{4\pi}\,{}_2F_1(-2k\!-\!1,1;\f{3}{2};2)
\ee
for $k=0, 1, 2, \ldots\ $. In addition, they gave the approximation \cite[(1.26)]{DRZ}
\[4\alpha^{3/4}\int_0^\infty\frac{xe^{-\alpha x^2}}{e^{2\pi x}-1}\,{}_1F_1(-2k;\f{3}{2};2\alpha x^2)\,dx\hspace{4cm}\]
\bee\label{e13}
\simeq {}_2F_1(-2k,1;\f{3}{2};2) \bl\{-\alpha^{-1/4}+\bl(\frac{1}{\alpha}+\frac{1}{\beta}+\frac{2}{3\cdot{}_2F_1(-2k,1;\f{3}{2};2)}\br)^{1/4}\br\},
\ee
which reduces to (\ref{e11}) when $k=0$.

At the end of their paper, Dixit {\it et al.\/} posed the following two open questions, namely:

{\bf Question 1.}\ Find the exact evaluation of the integral
\bee\label{e15}
{\hat J}_{2k}(\pi):=\int_0^\infty\frac{xe^{-\pi x^2}}{e^{2\pi x}-1}\,{}_1F_1(-2k;\f{3}{2};2\pi x^2)\,dx
\ee
for positive integer $k$.
\vspace{0.2cm}

{\bf Question 2.}\  Find the exact evaluation of, or at least an approximation to, the integral
\bee\label{e16}
{\hat J}_{2k+1}(\alpha)=\int_0^\infty\frac{xe^{-\alpha x^2}}{e^{2\pi x}-1}\,{}_1F_1(-2k\!-\!1;\f{3}{2};2\alpha x^2)\,dx 
\ee
when $\alpha\neq\pi$ is a positive real number and $k$ is a non-negative integer.
In this note we partially answer the above two questions by obtaining simple closed-form expressions for these integrals which, although not exact, approximate the given integrals to within exponentially small accuracy when $k$ is large and $a=O(1)$. In addition, we extend the scope of Question 1 by considering the integral ${\hat J}_{2k}(\alpha)$ with $\alpha>0$ and, as a by-product of the analysis pertaining to Question 2, we
supply an alternative proof of the result (\ref{e14}). 
\vspace{0.6cm}

\begin{center}
{\bf 2.\ The analysis of $J_{2k}(a)$} 
\end{center}
\setcounter{section}{2}
\setcounter{equation}{0}
\renewcommand{\theequation}{\arabic{section}.\arabic{equation}}
Throughout we shall find it convenient to replace the parameter $\alpha$ by $\pi a$ and define
the integral $J_{2k}(a)$  by
\bee\label{e20}
J_{2k}(a)=\int_0^\infty\frac{xe^{-\pi ax^2}}{e^{2\pi x}-1}\,{}_1F_1(-2k;\f{3}{2};2\pi ax^2)\,dx
\ee
for $a>0$ and positive integer $k$. Then $J_{2k}(1)={\hat J}_{2k}(\pi)$ in (\ref{e15}).
The confluent hypergeometric function terminates and we have \cite[p. 322]{DLMF}
\[{}_1F_1(-2k;\f{3}{2};2\pi x^2)=\sum_{r=0}^{2k}\frac{(-2k)_r}{(\f{3}{2})_r r!}(2\pi ax^2)^r.\]
Substitution of this series into the left-hand side of the above yields
\[J_{2k}(a)
=\sum_{r=0}^{2k}\frac{(-2k)_r}{(\f{3}{2})_r r!} (2\pi a)^r \int_0^\infty\frac{x^{2r+1}e^{-\pi ax^2}}{e^{2\pi x}-1}\,dx\]
upon reversal of the order of summation and integration. 

Now
\[\int_0^\infty\frac{x^{2r+1}e^{-\pi ax^2}}{e^{2\pi x}-1}\,dx=\sum_{n\geq 1} \int_0^\infty x^{2r+1} e^{-\pi ax^2-2\pi nx}dx=\frac{a^{1/2}r! \g(r+\f{3}{2})}{2(\pi a)^{r+3/2}}\,\,U_r,\]
where 
\[U_r:=\sum_{n\geq 1} U(r+1,\f{1}{2},\pi n^2/a)\]
with $U(a,b,z)$ being the confluent hypergeometric function of the second kind \cite[p. 322]{DLMF}. Then we obtain 
\[J_{2k}(a)=\frac{1}{4\pi a}\sum_{r=0}^{2k} (-2k)_r 2^rU_r.\]

From the integral representation \cite[p. 326]{DLMF}
\bee\label{e21b}
U(a,b,z)=\frac{1}{\g(a)}\int_0^\infty e^{-zt} t^{a-1}(1+t)^{b-a-1}dt\qquad (a>0,\,\Re (z)>0),
\ee
we find 
\[U_r=\sum_{n\geq 1} \frac{1}{r!} \int_0^\infty e^{-\pi n^2t/a} t^r(1+t)^{-r-3/2}dt=\frac{1}{r!}\int_0^\infty\!\!\!
\frac{\psi(t)\,t^r}{(1+t)^{r+3/2}}\,dt,\]
where we have defined
\bee\label{e21a}
\psi(t):=\sum_{n\geq 1} e^{-\pi n^2t/a}.
\ee
Hence
\[J_{2k}(a)=\frac{1}{4\pi a}\int_0^\infty \sum_{r=0}^{2k}\frac{(-2k)_r}{r!} \bl(\frac{2t}{1+t}\br)^r \frac{\psi(t)}{(1+t)^{3/2}}\,dt\]
\bee\label{e21}
=\frac{1}{4\pi a}\int_0^\infty \frac{\psi(t) (1-t)^{2k}}{(1+t)^{2k+3/2}}\,dt,\hspace{2cm}
\ee
where the finite sum has been evaluated as \cite[(15.4.6)]{DLMF}
\[{}_1F_0\bl(-2k;;\frac{2t}{1+t}\br)=\bl(\frac{1-t}{1+t}\br)^{2k}.\]

We now divide the integration path into $[0,1]$ and $[1,\infty)$ and make the change of variable $t\to 1/t$ in the integral over $[0,1]$. This yields
\[J_{2k}(a)=\frac{1}{4\pi a}\int_1^\infty \{t^{-1/2} \psi(1/t)+\psi(t)\}\,\frac{(t-1)^{2k}}{(1+t)^{2k+3/2}}\,dt.\]
For the sum 
\bee\label{e21c}
\Psi(\tau)=\sum_{n\geq 1} e^{-\pi n^2\tau},
\ee
we have the well-known Poisson transformation given by \cite[p. 124]{WW}
\bee\label{e21d}
\Psi(\tau)+\fs(1-\tau^{-1/2})=\tau^{-1/2}\Psi(1/\tau).
\ee
With $\tau=at$, this yields
\bee\label{e22}
t^{-1/2} \psi(1/t)=a^{1/2}\{\phi(t)+\fs(1-(at)^{-1/2})\}, \qquad \phi(t):=\sum_{n\geq 1} e^{-\pi n^2at}.
\ee
Hence
\begin{eqnarray*}
J_{2k}(a)&=&\frac{1}{4\pi a}\int_1^\infty\bl\{\psi(t)+a^{1/2}\phi(t)+\fs a^{1/2}(1-(at)^{-1/2})\br\}\,\frac{(t-1)^{2k}}{(1+t)^{2k+3/2}}\,dt\\
&=&-\frac{1}{8\pi a}\int_0^1(1-a^{1/2}t^{-1/2})\,\frac{(1-t)^{2k}}{(1+t)^{2k+3/2}}\,dt\\
&&\hspace{4cm}+\frac{1}{4\pi a}\int_1^\infty\bl\{\psi(t)+a^{1/2}\phi(t)\br\}\,\frac{(t-1)^{2k}}{(1+t)^{2k+3/2}}\,dt.
\end{eqnarray*}

For positive integer $k$, we have the integrals
\[\int_0^1 \frac{t^{-1/2}(1-t)^{2k}}{(1+t)^{2k+3/2}}\,dt=\sqrt{\frac{\pi}{2}}\,\frac{\g(2k+1)}{\g(2k+\f{3}{2})}\]
and
\begin{eqnarray*}
\int_0^1\frac{(1-t)^{2k}}{(1+t)^{2k+3/2}}\,dt&=&\frac{1}{2k+1}\,{}_2F_1(1,2k+\f{3}{2};2k+2;-1)\\
&=&2\,{}_2F_1(-2k,1;\f{3}{2};2)-\sqrt{\frac{\pi}{2}}\,\frac{\g(2k+1)}{\g(2k+\f{3}{2})}
\end{eqnarray*}
by application of the transformation \cite[p. 390]{DLMF}
\[{}_2F_1(a,b;c;z)=\frac{\g(a)\g(c-a-b)}{\g(c-a)\g(c-b)}\,z^{-a} {}_2F_1(a,a-c+1;a+b-c+1;1-z^{-1})\]
\bee\label{e25}
+\frac{\g(c)\g(a+b-c)}{\g(a)\g(b)}\,z^{a-c}(1-z)^{c-a-b} {}_2F_1(c-a,1-a;c-a-b+1;1-z^{-1}).
\ee
Hence we obtain
\newtheorem{theorem}{Theorem}
\begin{theorem}$\!\!\!.$\ Let $a>0$ and $k$ be a positive integer. Then the integral $J_{2k}(a)$ defined in (\ref{e20}) satisfies
\bee\label{e26}
J_{2k}(a)=T_{2k}(a)+\epsilon_{2k}(a),
\ee
where
\bee\label{e27a}
T_{2k}(a)=\frac{1}{4\pi a}\bl\{\bl(\frac{1+a^{1/2}}{2}\br)\sqrt{\frac{\pi}{2}}\,\frac{\g(2k+1)}{\g(2k+\f{3}{2})}-{}_2F_1(-2k,1;\f{3}{2};2)
\br\}
\ee
and
\bee\label{e27}
\epsilon_{2k}(a)=\frac{1}{4\pi a}\int_1^\infty \{\psi(t)+a^{1/2}\phi(t)\}\frac{(t-1)^{2k}}{(1+t)^{2k+3/2}}\,dt
\ee
with the sums $\psi(t)$ and $\phi(t)$ defined in (\ref{e21a}) and (\ref{e22}).
\end{theorem}

It will be found subsequently that $\epsilon_{2k}(a)$ is small for $k\geq 1$ when $a=O(1)$ and so we shall refer to it as the remainder term. We observe that when $a=1$, we have $\phi(t)=\psi(t)$ and hence that
\[\epsilon_{2k}(1)=\frac{1}{2\pi}\int_1^\infty \frac{\psi(t) (t-1)^{2k}}{(1+t)^{2k+3/2}}\,dt.\]
\vspace{0.6cm}

\begin{center}
{\bf 3.\ The analysis of $J_{2k+1}(a)$} 
\end{center}
\setcounter{section}{3}
\setcounter{equation}{0}
\renewcommand{\theequation}{\arabic{section}.\arabic{equation}}
A similar treatment for the integral
\bee\label{e31}
J_{2k+1}(a)=\int_0^\infty\frac{xe^{-\pi ax^2}}{e^{2\pi x}-1}\,{}_1F_1(-2k\!-\!1;\f{3}{2};2\pi ax^2)\,dx\qquad (k=0, 1, 2, \ldots)
\ee
shows that
\[J_{2k+1}(a)=\frac{1}{4\pi a}\sum_{r=0}^{2k+1}(-2k\!-\!1)_r 2^r U_r=\frac{1}{4\pi a}\int_0^\infty \frac{\psi(t)\,(1-t)^{2k+1}}{(1+t)^{2k+5/2}}\,dt.\]
Dividing the integration path as in Section 2, we find 
\[J_{2k+1}(a)=\frac{1}{4\pi a}\int_1^\infty\{t^{-1/2}\psi(1/t)-\psi(t)\}\,\frac{(t-1)^{2k+1}}{(1+t)^{2k+5/2}}\,dt.\]
Application of (\ref{e22}) and some straightforward algebra then produces
\[J_{2k+1}(a)=-\frac{1}{8\pi a}\int_0^1(1-a^{1/2}t^{-1/2})\,\frac{(1-t)^{2k+1}}{(1+t)^{2k+5/2}}\,dt\hspace{4cm}\]
\[\hspace{4cm}+\frac{1}{4\pi a}\int_1^\infty
\{a^{1/2}\phi(t)-\psi(t)\}\,\frac{(t-1)^{2k+1}}{(1+t)^{2k+5/2}}\,dt.\]

Now
\[\int_0^1(1-a^{1/2}t^{-1/2})\,\frac{(1-t)^{2k+1}}{(1+t)^{2k+5/2}}\,dt=\frac{1}{2k+2}\,{}_2F_1(1,2k+\f{5}{2};2k+3;-1)-\sqrt{\frac{\pi a}{2}}\,\frac{\g(2k+2)}{\g(2k+\f{5}{2})}\]
\[=2\,{}_2F_1(-2k\!-\!1,1;\f{3}{2};2)+(1-a^{1/2})\sqrt{\frac{\pi}{2}}\,\frac{\g(2k+2)}{\g(2k+\f{5}{2})}\]
by (\ref{e25}). Hence we obtain
\begin{theorem}$\!\!\!.$\ Let $a>0$ and $k$ be a non-negative integer. Then the integral $J_{2k+1}(a)$ defined in (\ref{e31}) satisfies
\bee\label{e32a}
J_{2k+1}(a)=-T_{2k+1}(a)+\epsilon_{2k+1}(a),
\ee
where
\bee\label{e32c}
T_{2k+1}(a)=\frac{1}{4\pi a}\bl\{\bl(\frac{1-a^{1/2}}{2}\br) \sqrt{\frac{\pi}{2}}\,\frac{\g(2k+2)}{\g(2k+\f{5}{2})}+{}_2F_1(-2k\!-\!1,1;\f{3}{2};2)\br\}
\ee
and
\bee\label{e32b}
\epsilon_{2k+1}(a)=\frac{1}{4\pi a}\int_1^\infty \{a^{1/2}\phi(t)-\psi(t)\}\frac{(t-1)^{2k+1}}{(1+t)^{2k+5/2}}\,dt
\ee
with the sums $\psi(t)$ and $\phi(t)$ defined in (\ref{e21a}) and (\ref{e22}).
\end{theorem}

When $a=1$, we have $\psi(t)=\phi(t)$ and hence $\epsilon_{2k+1}(1)=0$. It then follows from (\ref{e32a}) that
\[J_{2k+1}(1)=-\frac{1}{4\pi}\,{}_2F_1(-2k\!-\!1,1;\f{3}{2};2),\]
which supplies another proof of the result stated in (\ref{e14}) obtained in \cite{DRZ}.

\vspace{0.6cm}

\begin{center}
{\bf 4.\ Estimation of the remainder terms} 
\end{center}
\setcounter{section}{4}
\setcounter{equation}{0}
\renewcommand{\theequation}{\arabic{section}.\arabic{equation}}
We examine the remainder terms $\epsilon_{2k}(a)$ and $\epsilon_{2k+1}(a)$ appearing in (\ref{e27}) and (\ref{e32b})
and determine bounds and an estimate of their large-$k$ behaviour. We consider first the term $\epsilon_{2k}(a)$
which can be written as
\[\epsilon_{2k}(a)=\frac{a^{-3/4}}{4\pi}\int_1^\infty\{a^{1/4}\phi(t)+a^{-1/4}\psi(t)\}\,\frac{(t-1)^{2k}}{(1+t)^{2k+3/2}}\,dt\]
With the change of variable $t\to 1+u$, we have
\begin{eqnarray}
\epsilon_{2k}(a)\!\!&=&\!\!\frac{a^{-3/4}}{4\pi}\bl\{a^{1/4}\sum_{n\geq 1}e^{-\pi n^2a}\!\!\int_0^\infty \!\!\!e^{-\pi n^2au}h(u)\,du+
a^{-1/4}\sum_{n\geq 1}e^{-\pi n^2/a}\!\!\int_0^\infty \!\!\!e^{-\pi n^2u/a}h(u)\,du\br\}\nonumber\\
&<&\frac{a^{-3/4}}{4\pi}\bl\{a^{1/4}\Psi(a)\int_0^\infty \!\!\!e^{-\pi au} h(u)\,du+a^{-1/4}\Psi(1/a)\int_0^\infty\!\!\!e^{-\pi u/a}h(u)\,du\br\},\label{e41}
\end{eqnarray}
where $\Psi(a)$ is defined in (\ref{e21c}) and
$h(u)=u^{2k}/(2+u)^{2k+3/2}$.
Evaluation of the integrals appearing in (\ref{e41}) in terms of the confluent hypergeometric function $U(a,b,z)$ by (\ref{e21b}), we then obtain the upper bound in the form
\begin{theorem}$\!\!\!.$\ The remainder term $\epsilon_{2k}(a)$ defined in (\ref{e27}) satisfies the upper bound
\bee\label{e42}
\epsilon_{2k}(a)<{\cal B}_{2k}(a),\qquad {\cal B}_{2k}(a):=\frac{a^{-3/4}(2k)!}{4\sqrt{2} \pi}\{E_{2k}(a)+E_{2k}(1/a)\},
\ee
where
\[E_{2k}(a):=a^{1/4}\Psi(a) U(2k+1,\fs,2\pi a).\]
and $\Psi(a)$ is given by (\ref{e21c}). 
\end{theorem}

The behaviour of this bound as $k\to\infty$ with $a$ fixed can be obtained by making use of the result \cite[(13.8.8)]{DLMF}
\[U(2k+1,\fs,2\pi a)\sim \frac{e^{\pi a}}{(2k)!} \sqrt{\frac{\pi}{2k}}\,e^{-4\sqrt{\pi ak}}\qquad (k\to\infty,\ a\ll 2k/\pi).\]
For values of $a\simeq 1$, we can bound the sum $\Psi(a)$ by 
\[\Psi(a):=\sum_{n\geq 1} e^{-\pi n^2a}=e^{-\pi a}\bl(1+e^{-3\pi a}+e^{\pi a} \sum_{n\geq 3}e^{-\pi n^2a}\br)<\lambda(a) e^{-\pi a},\]
where
\bee\label{e33a}
\lambda(a):=1+e^{-3\pi a}+e^{\pi a} \sum_{n\geq 3}e^{-\pi na}=1+e^{-3\pi a}+\frac{e^{-2\pi a}}{1-e^{-\pi a}}.
\ee
This then yields the estimate as $k\to\infty$
\bee\label{e43}
{\cal B}_{2k}(a)\sim \frac{a^{-3/4}k^{-1/2}}{8\sqrt{\pi}}\bl\{a^{1/4}\lambda(a) e^{-4\sqrt{\pi ak}}+a^{-1/4}\lambda(1/a) e^{-4\sqrt{\pi k/a}}\br\}
\ee
provided $a\gg \pi/(2k)$ and $a\ll 2k/\pi$ (that is, when $a$ is neither too small nor too large).  In the case $a=1$ we have
\[{\cal B}_{2k}(1)\sim\frac{\lambda(1)}{4\sqrt{\pi}}\,k^{-1/2} e^{-4\sqrt{\pi k}}\qquad (k\to\infty).\]

The remainder term $\epsilon_{2k+1}(a)$ may be written as
\[\epsilon_{2k+1}(a)=\frac{a^{-3/4}}{4\pi}\int_1^\infty \{a^{1/4}\phi(t)-a^{-1/4}\psi(t)\}\,\frac{(t-1)^{2k+1}}{(1+t)^{2k+5/2}}\,dt.\]
It is straightforward to show (we omit these details) that $a^{1/2}\phi(t)-\psi(t)$ has opposite signs in the intervals $a\in(0,1)$ and $a\in(1,\infty)$ when $t\in[1,\infty)$, being negative in $a\in(1,\infty)$. Hence it follows that $\epsilon_{2k+1}(a)<0$ when $a\in(1,\infty)$ and $\epsilon_{2k+1}(a)>0$ when $a\in(0,1)$. The same procedure employed for $\epsilon_{2k}(a)$ shows that\footnote{It is clear that this bound will not be sharp in the neighbourhood of $a\simeq 1$.}
\[|\epsilon_{2k+1}(a)|<\frac{a^{-3/4}}{4\pi}\int_1^\infty \{a^{1/4}\phi(t)+a^{-1/4}\psi(t)\}\,\frac{(t-1)^{2k+1}}{(1+t)^{2k+5/2}}\,dt\]
and therefore we obtain
\begin{theorem}$\!\!\!.$\ The remainder term $\epsilon_{2k+1}(a)$ defined in (\ref{e32b}) satisfies the upper bound
\bee\label{e44}
|\epsilon_{2k+1}(a)|<{\cal B}_{2k+1}(a),\qquad {\cal B}_{2k+1}(a):=\frac{a^{-3/4}(2k+1)!}{4\sqrt{2} \pi}\{E_{2k+1}(a)+E_{2k+1}(1/a)\},
\ee
where
\[E_{2k+1}(a):=a^{1/4}\Psi(a) U(2k+2,\fs,2\pi a).\]
and $\Psi(a)$ is given by (\ref{e21c}). The leading behaviour of ${\cal B}_{2k+1}(a)$ for large $k$ and finite $a$ is given by the right-hand side of (\ref{e43}).
\end{theorem}  
\bigskip
\vspace{0.6cm}

\begin{center}
{\bf 5.\ Numerical results} 
\end{center}
\setcounter{section}{5}
\setcounter{equation}{0}
\renewcommand{\theequation}{\arabic{section}.\arabic{equation}}

To demonstrate the smallness of the remainder terms $\epsilon_{2k}(a)$ and $\epsilon_{2k+1}(a)$ we define the quantities
\[{\cal J}_{2k}(a):=J_{2k}(a)-\frac{1}{4\pi a}\bl\{\bl(\frac{1+a^{1/2}}{2}\br)     
\sqrt{\frac{\pi}{2}}\,\frac{\g(2k+1)}{\g(2k+\f{3}{2})}-{}_2F_1(-2k,1;\f{3}{2};2)\br\}\]
and
\[{\cal J}_{2k+1}(a):=J_{2k+1}(a)+\frac{1}{4\pi a}\bl\{\bl(\frac{1-a^{1/2}}{2}\br)\sqrt{\frac{\pi}{2}}\,\frac{\g(2k+2)}{\g(2k+\f{5}{2})}+{}_2F_1(-2k-1,1;\f{3}{2};2)\br\}.\]
In Tables 1--3 we present numerical values of these quantities compared with their bounds ${\cal B}_{2k}(a)$ and ${\cal B}_{2k+1}(a)$
for a range of $k$ and three values of the parameter $a=O(1)$. 
It is seen that this bound agrees very well with the computed values of ${\cal J}_{2k}(a)$ and ${\cal J}_{2k+1}(a)$. The estimates in (\ref{e43}) and (\ref{e44}) show that the remainder terms are {\it exponentially small} for large $k$ when $a=O(1)$. Consequently, the terms $T_{2k}(a)$ and $T_{2k+1}(a)$ in (\ref{e27a}) and (\ref{e32c}) approximate $J_{2k}(a)$ and $J_{2k+1}(a)$, respectively, to exponential accuracy in the large-$k$ limit. 

Now
\[J_{2k}(0)=J_{2k+1}(0)=\int_0^\infty\frac{x}{e^{2\pi x}-1}\,dx=\frac{1}{24};\]
but it is easily seen that $T_{2k}(a)$ and $T_{2k+1}(a)$ in (\ref{e27a}) and (\ref{e32c}) are $O(a^{-1})$ as $a\to 0$ and $O(a^{-1/2})$ as $a\to\infty$. Routine calculations show that the bounds ${\cal B}_{2k}(a)$ and ${\cal B}_{2k+1}(a)$ also possess the same behaviour in these limits. Consequently, the approximtions $T_{2k}(a)$ and $T_{2k+1}(a)$ will not be good for small or large values of the parameter $a$, although it is worth pointing out that the range of validity in $a$ will increase as $k$ increases.

\begin{table}[h]
\caption{\footnotesize{Values of ${\cal J}_{2k}(a)$ and the bound ${\cal B}_{2k}(a)$ for $\epsilon_{2k}(a)$ in (\ref{e42}) as a function of $k$ when $a=1$.}}
\begin{center}
\begin{tabular}{l|c|c||l|c|c}
\mcol{1}{c|}{$k$} & \mcol{1}{c|}{${\cal J}_{2k}(1)$} & \mcol{1}{c||}{${\cal B}_{2k}(1)$ } & \mcol{1}{c|}{$k$}
& \mcol{1}{c|}{${\cal J}_{2k}(1)$} & \mcol{1}{c}{${\cal B}_{2k}(1)$}\\
[0.05cm]\hline
&&&&&\\[-0.25cm]
1 & $1.250\times 10^{-5}$ & $1.253\times 10^{-5}$ & 10 & $3.905\times 10^{-12}$ & $3.913\times 10^{-12}$\\
2 & $8.571\times 10^{-7}$ & $8.588\times 10^{-7}$ & 20 & $3.186\times 10^{-16}$ & $3.193\times 10^{-16}$\\
3 & $9.818\times 10^{-8}$ & $9.838\times 10^{-8}$ & 30 & $2.305\times 10^{-19}$ & $2.309\times 10^{-19}$\\
5 & $2.883\times 10^{-9}$ & $2.888\times 10^{-9}$ & 50 & $2.433\times 10^{-24}$ & $2.438\times 10^{-24}$\\
[.10cm]\hline
\end{tabular}
\end{center}
\end{table}

\begin{table}[h]
\caption{\footnotesize{Values of ${\cal J}_{2k}(a)$ and the bound ${\cal B}_{2k}(a)$ as a function of $k$ when $a=2$ and $a=0.50$.}}
\begin{center}
\begin{tabular}{l|l|l||l|l|l}
\mcol{1}{c|}{$k$} & \mcol{1}{c|}{${\cal J}_{2k}(2)$} & \mcol{1}{c||}{${\cal B}_{2k}(2)$ } & \mcol{1}{c|}{$k$}
& \mcol{1}{c|}{${\cal J}_{2k}(2)$} & \mcol{1}{c}{${\cal B}_{2k}(2)$}\\
[0.05cm]\hline
&&&&&\\[-0.25cm]
1 & $5.987\times 10^{-5}$ & $6.364\times 10^{-5}$ & 10 & $9.509\times 10^{-10}$ & $1.011\times 10^{-9}$\\
2 & $7.856\times 10^{-6}$ & $8.355\times 10^{-6}$ & 20 & $1.075\times 10^{-12}$ & $1.143\times 10^{-12}$\\
3 & $1.563\times 10^{-6}$ & $1.662\times 10^{-6}$ & 30 & $6.016\times 10^{-15}$ & $6.398\times 10^{-15}$\\
5 & $1.162\times 10^{-7}$ & $1.236\times 10^{-7}$ & 50 & $1.668\times 10^{-18}$ & $1.774\times 10^{-18}$\\
[.10cm]\hline
\mcol{1}{c|}{$k$} & \mcol{1}{c|}{${\cal J}_{2k}(\fs)$} & \mcol{1}{c||}{${\cal B}_{2k}(\fs)$ } & \mcol{1}{c|}{$k$}
& \mcol{1}{c|}{${\cal J}_{2k}(\fs)$} & \mcol{1}{c}{${\cal B}_{2k}(\fs)$}\\
[0.05cm]\hline
&&&&&\\[-0.25cm]
1 & $1.693\times 10^{-4}$ & $1.800\times 10^{-4}$ & 10 & $2.689\times 10^{-9}$ & $2.860\times 10^{-9}$\\
2 & $2.222\times 10^{-5}$ & $2.363\times 10^{-5}$ & 20 & $3.040\times 10^{-12}$ & $3.234\times 10^{-12}$\\
3 & $4.420\times 10^{-6}$ & $4.700\times 10^{-6}$ & 30 & $1.702\times 10^{-14}$ & $1.810\times 10^{-14}$\\
5 & $3.287\times 10^{-7}$ & $3.496\times 10^{-7}$ & 50 & $4.719\times 10^{-18}$ & $5.019\times 10^{-18}$\\
[.10cm]\hline
\end{tabular}
\end{center}
\end{table}
\begin{table}[h]
\caption{\footnotesize{Values of ${\cal J}_{2k+1}(a)$ and the bound ${\cal B}_{2k+1}(a)$ as a function of $k$ when $a=2$ and $a=0.50$.}}
\begin{center}
\begin{tabular}{l|l|l||l|l|l}
\mcol{1}{c|}{$k$} & \mcol{1}{c|}{${\cal J}_{2k+1}(2)$} & \mcol{1}{c||}{${\cal B}_{2k+1}(2)$ } & \mcol{1}{c|}{$k$}
& \mcol{1}{c|}{${\cal J}_{2k+1}(2)$} & \mcol{1}{c}{${\cal B}_{2k+1}(2)$}\\
[0.05cm]\hline
&&&&&\\[-0.25cm]
0 & $2.230\times 10^{-4}$ & $2.376\times 10^{-4}$ & 10 & $6.340\times 10^{-10}$ & $6.743\times 10^{-10}$\\
1 & $2.018\times 10^{-5}$ & $2.147\times 10^{-5}$ & 20 & $8.067\times 10^{-13}$ & $8.579\times 10^{-13}$\\
2 & $3.376\times 10^{-6}$ & $3.591\times 10^{-6}$ & 30 & $4.760\times 10^{-15}$ & $5.063\times 10^{-15}$\\
5 & $6.603\times 10^{-8}$ & $7.022\times 10^{-8}$ & 40 & $6.287\times 10^{-17}$ & $6.686\times 10^{-17}$\\
[.10cm]\hline
\mcol{1}{c|}{$k$} & \mcol{1}{c|}{${\cal J}_{2k+1}(\fs)$} & \mcol{1}{c||}{${\cal B}_{2k+1}(\fs)$ } & \mcol{1}{c|}{$k$}
& \mcol{1}{c|}{${\cal J}_{2k+1}(\fs)$} & \mcol{1}{c}{${\cal B}_{2k+1}(\fs)$}\\
[0.05cm]\hline
&&&&&\\[-0.25cm]
0 & $6.307\times 10^{-4}$ & $6.720\times 10^{-4}$ & 10 & $1.793\times 10^{-9}$ & $1.907\times 10^{-9}$\\
1 & $5.708\times 10^{-5}$ & $6.072\times 10^{-5}$ & 20 & $2.282\times 10^{-12}$ & $2.427\times 10^{-12}$\\
2 & $9.548\times 10^{-6}$ & $1.106\times 10^{-5}$ & 30 & $1.346\times 10^{-14}$ & $1.432\times 10^{-14}$\\
5 & $1.867\times 10^{-7}$ & $1.986\times 10^{-7}$ & 40 & $1.778\times 10^{-16}$ & $1.891\times 10^{-16}$\\
[.10cm]\hline
\end{tabular}
\end{center}
\end{table}

The approximation in (\ref{e13}) when $\alpha=a\pi$ (with $\alpha\beta=\pi^2$) yields
\bee\label{e51}
J_{2k}(a)\simeq -\frac{F}{4\pi a}\,\bl\{1-\bl(1+a^2+\frac{2\pi a}{3F}\br)^{1/4}\br\},\qquad F:={}_2F_1(-2k,1;\f{3}{2};2).
\ee
This yields the limiting behaviours
\[J_{2k}(a)\simeq\frac{1}{24}+a\bl(\frac{F}{16\pi}-\frac{\pi}{96F}\br)+O(a^2) \qquad (a\to 0)\]
and 
\[J_{2k}(a)\simeq \frac{F}{4\pi\sqrt{a}}\,\bl\{1-\frac{1}{\sqrt{a}}+\frac{\pi}{6aF}+O(a^{-3/2})\br\} \qquad (a\to\infty).\]
The approximation (\ref{e51}) is found to be quite accurate in the limits of small and large $a$, with $k$ finite. However, the accuracy is not good when $a=O(1)$. For example, when $a=1$ and
$k=5, 10$ the approximation (\ref{e51}) yields absolute relative errors of 8.8\% and 19.2\%, respectively; and this error increases as $k$ increases, in marked contrast to the approximations in (\ref{e27a}) and (\ref{e32c}).

As a final remark, it is doubtful that the remainder terms $\epsilon_{2k}(a)$ and $\epsilon_{2k+1}(a)$ can be expressed in simple closed forms.

\end{document}